\documentclass[12pt]{article}
\usepackage{amssymb}
\usepackage[mathscr]{eucal}
\usepackage{amsthm}
\usepackage{eufrak}
\usepackage{amstext}
\usepackage{amsmath}

\begin{document}

\def\cC{{\cal C}}
\def\cF{{\cal F}}
\def\cH{{\cal H}}
\def\cP{{\cal P}}
\def\bd{\mathop{\rm bd}}
\def\cl{\mathop{\rm cl}}
\def\rdiv{\mathop{\rm div}}
\def\const{\mathop{\rm const}}
\def\mn{{\medskip\noindent}}

\newtheorem{theorem}{Theorem}{\bf}{\it}
\newenvironment{form}{\begin{list}{}{\labelwidth 1.1cm
\labelsep 0.2cm \leftmargin 1.3cm}}{\end{list}}

\begin{center}
{\Large\bf Lattice Points in Large Borel Sets and Successive
Minima}
\end{center}

\begin{center}
Iskander Aliev and Peter M.~Gruber
\end{center}

\begin{quote}
\footnotesize{ \vskip 1truecm\noindent {\bf Abstract.} Let $B$ be
a Borel set in $\mathbb E^{d}$ with volume $V(B)=\infty$. It is
shown that almost all lattices $L$ in $\mathbb E^{d}$ contain
infinitely many pairwise disjoint $d$-tuples, that is sets of $d$
linearly independent points in $B$. A consequence of this result
is the following: let $S$ be a star body in $\mathbb E^{d}$ with
$V(S)=\infty$. Then for almost all lattices $L$ in $\mathbb
E^{d}$ the successive minima $\lambda_{1}(S,L),\dots ,
\lambda_{d}(S,L)$ of $S$ with respect to $L$ are $0$. A
corresponding result holds for most lattices in the Baire
category sense. A tool for the latter result is the
semi-continuity of the successive minima.


}

\end{quote}

\section{Introduction and Statement of Results}

A {\it lattice} $L$ in Euclidean $d$-space $\mathbb E^{d}$ is the
system of all integer linear combinations of $d$ linearly
independent vectors in $\mathbb E^{d}$. These vectors form a {\it
basis} of $L$ and the absolute value of their determinant is the
{\it determinant} $d(L)$ of $L$. $d(L)$ is independent of the
particular choice of a basis. To each lattice $L$ we let
correspond all $d\times d$ matrices the column vectors of which
form a basis of $L$. Identify each such matrix with a point in
$\mathbb E^{d^2}$. There is a Borel set $\mathscr F$ in $\mathbb
E^{d^2}$ consisting of such matrices, which has infinite Lebesgue
measure and such that to each lattice $L$ corresponds precisely
one matrix in $\mathscr F$. Thus there is a one-to-one
correspondence between the space $\mathscr L$ of all lattices in
$\mathbb E^d$ and $\mathscr F$. The Lebesgue measure on $\mathscr
F$ then yields a measure $\nu$ on $\mathscr L$.

Results of Rogers \cite{Rogers55} (for $d\ge 3)$ and Schmidt
\cite{Schmidt60} (for $d=2)$ show that for a Borel set $B$ in
$\mathbb E^{d}$ with Lebesgue measure $V(B)=\infty$, $\nu$-almost
all lattices $L$ contain infinitely many primitive points in $B$,
where a point $l\in L$ is {\it primitive} if it is different from
the origin $o$ and on the line-segment $[o,l]$ there are no
points of $L$, except $o$ and $l$. A refinement of this result is
as follows.

\begin{theorem}
Let $B$ be a Borel set in $\mathbb E^{d}$ with $V(B)=\infty$.
Then for $\nu$-almost every lattice $L\in\mathscr L$, the set $B$
contains infinitely many, pairwise disjoint $d$-tuples of
linearly independent primitive points of $L$.
\end{theorem}
\noindent Tools for the proof are measure theoretic results of
Rogers \cite{Rogers55} and Schmidt \cite{Schmidt60} and a result
of Yao and Yao \cite{YaoYao85} from applied computational
geometry on dissection of sets in $\mathbb E^{d}$.

A {\it star body} $S$ in $\mathbb E^{d}$ is a closed set with $o$
in its interior such that each ray with endpoint $o$ meets the
boundary of $S$ in at most one point. Equivalently, $S=\{
x:f(x)\le 1\}$, where $f:\mathbb E^{d}\to\mathbb R$ is a {\it
distance function}, i.e. it is non-negative, continuous and
positively homogeneous of degree $1$. The {\it successive minima}
of $S$ or $f$ with respect to a lattice $L$ are defined as
follows:
\begin{align*}
\lambda_{i}(S,L)&=\lambda_{i}(f,L)\\
&=\inf\{\lambda >0:\lambda S\cap L \text{ contains } i \text{
linearly independent vectors}\}\\
&=\inf\{\max\{f(l_{1}),\dots , f(l_{i})\}:l_{1},\dots , l_{i}\in L
\text{ linearly independent}\}
\end{align*}
for $i=1,\dots,d$. Clearly,

\begin{form}
\item[\hbox to 1.1cm{\rm (1)\hfill}]
$0\le\lambda_{1}(S,L)\le\cdots\le\lambda_{d}(S,L)\le\infty .$
\end{form}

Successive minima play an important role in the geometry of
numbers, algebraic number theory, Diophantine approximation and
computational geometry, see e.g.\
\cite{GruberLekkerkerker87,Gruber93,Schmidt80,Lagarias95,Blomer00}.
For a surprising relation to Nevanlinna's value distribution
theory see \cite{Hyuga00}.

Let $\mathscr L$ be endowed with its natural topology, see
\cite{GruberLekkerkerker87}. Then $\mathscr L$ is locally compact
by Mahler's compactness theorem. Thus a version of the Baire
category theorem implies that $\mathscr L$ is {\it Baire}. That
is, any meager set has dense complement, where a set is {\it
meager} or {\it of first Baire category}, it it is a countable
union of nowhere dense sets, see \cite{Oxtoby80}.

\begin{theorem}
Let $S$ be a star body in $\mathbb E^{d}$ with $V(S)=\infty$. Then
$\lambda_{1}(S,L)=\cdots \\=\lambda_{d}(S,L)=0$ for
\begin{itemize}
\item[{\rm (i)}] $\nu$-almost all lattices $L$ in $\mathscr L$ and for
\item[{\rm (ii)}] all lattices $L$ in $\mathscr L$, with a meager set
of exceptions.
\end{itemize}
\end{theorem}
\noindent Tools for the proof are Theorem~1 and a semi-continuity
result for successive minima which may be described as follows:

Let $(S_n)$ be a sequence of star bodies and $(f_n)$ the
corresponding sequence of distance functions. Then $(S_n)$ {\it
converges} to a star body $S$ with corresponding distance function
$f$ if the sequence $(f_n)$ converges uniformly to $f$ on the
solid unit ball $\{x:\Vert x\Vert \le 1\}$ of $\mathbb E^{d}$. A
sequence $(L_{n})$ of lattices {\it converges} to a lattice $L$,
if there are bases $\{ b_{n1},\dots , b_{nd}\}$ of $L_{n}$ and
$\{b_{1},\dots , b_{d}\}$ of $L$ such that $b_{n1}\to b_{1},\dots
, b_{nd}\to b_{d}$. This notion of convergence induces the
topology on $\mathscr L$.

\medskip\noindent{\bf Lemma.} {\it Let} $(S_{n})$ {\it be a sequence
of star bodies  and} $(L_n)$ {\it a sequence of lattices in}
$\mathbb E^{d}$, {\it converging to a star body} $S$ {\it and a
lattice $L$, respectively. Then,}

\begin{itemize}
\item[{\rm (i)}] $\displaystyle\limsup_{n\to\infty}
\lambda_{i}(S_{n},L_{n})\le\lambda_{i}(S,L),\,{\it for}\,\,
i=1,\dots , d,\,\,{\it and}$
\item[{\rm (ii)}] {\it if} $S$ {\it is bounded, then} $\displaystyle
\lim_{n\to\infty}\lambda_{i}(S_{n},L_{n})$ {\it exists and is
equal to} $\lambda_{i}(S,L),\\ {\it for}\,\,i=1,\dots ,d$.
\end{itemize}
\medskip

\noindent To see that $\lambda_{i}$ is {\it not} continuous, let
$S$ be a star body with $V(S)=\infty$ such that there is a
lattice $L$ which has only $o$ in common with the interior of
$S$, for example the star body $\{x:|x_{1}\cdots x_{d}|\le 1\}$,
see \cite{GruberLekkerkerker87}, p. 28. Then
$1\le\lambda_{i}(S,L)<\infty$, while by Theorem~2 there is a
sequence $(L_{n})$ of lattices such that $L_{n}\to L$ with
$\lambda_{i}(S,L_{n})=0$ for all $n$.

\section{Proof of Theorem~1}

A result of Yao and Yao \cite{YaoYao85} says that any mass
distribution in $\mathbb E^{d}$ with positive, continuous density
which tends rapidly to $0$ as $\Vert x\Vert\to\infty$, and of
total mass $V$, can be dissected into $2^{d}$ disjoint Borel
parts, each of mass $2^{-d}V$ and such that no hyperplane meets
all these $2^{d}$ masses. We need the following version of this
result:

\begin{form}
\item[\hbox to 1.1cm{\rm (2)\hfill}]
Let $A\subset\mathbb E^{d}$ be a bounded Borel set with volume
$V(A)>V>0$. Then $A$ contains $2^{d}$ pairwise disjoint Borel
subsets, each of volume $2^{-d}V$ and such that no
$(d-1)$-dimensional subspace of $\mathbb E^{d}$ meets each of
these $2^{d}$ sets.
\end{form}
To see (2), choose a compact set $C\subset A$ with $V(C)>V$. This
is possible by the inner regularity of Lebesgue measure. Next,
choose a continuous function $g:\mathbb E^{d}\to\mathbb R^{+}$
such that
$$g\ge \chi_{C}, \int\limits_{\mathbb E^{d}}
(g-\chi_{C})\,dx<2^{-d}(V(C)-V), \,\,g(x)\to 0 \text{ rapidly as }
\Vert x\Vert\to\infty ,$$ where $\chi_{C}$ is the characteristic
function of $C$. This is possible by the outer regularity of
Lebesgue measure and Urysohn's lemma. Let $F_{i}, i=1,\dots ,
2^{d},$ be a dissection of $\mathbb E^{d}$ for the density $g$ as
described by Yao and Yao such that
$$\int\limits_{F_{i}} g \,dx=2^{-d} \int\limits_{\mathbb E^{d}} g\,
dx\ge 2^{-d}V(C).$$ Then there is no $(d-1)$-dimensional subspace
of $\mathbb E^{d}$ which meets each of the sets $C\cap F_{i}$,
and for the respective volumes of these sets we have the following
estimate:
\begin{align*}
V(C\cap F_{i})&=\int\limits_{F_{i}} \chi_{C}
dx=\int\limits_{F_{i}} g\,
dx-\int\limits_{F_{i}}(g-\chi_{C})\,dx\ge \int\limits_{F_{i}} g\,
dx-\int\limits_{\mathbb E^{d}}
(g-\chi_{C})\,dx\\
&>2^{-d}V(C)-2^{-d}(V(C)-V)=2^{-d}V.
\end{align*}
This concludes the proof of (2).

For the proof of Theorem~2 assume first that $d\ge 3$. The
following result is an immediate consequence of a result of Rogers
\cite{Rogers55}, p. 286:

\begin{form}
\item[\hbox to 1.1cm{\rm (3)\hfill}]
Let $k=1,2,\dots ,$ and $A$ a Borel set in $\mathbb E^{d}$ with
$0<V(A)<\infty$. Then the function $\#^*(A\cap\,\cdot\, ):\mathscr
L\to \{ 0,1,\dots \}$, which counts the number of primitive
points of $L$ in $A$, is Borel measurable and
$$\int\limits_{\mathscr L(k)} \Big(\#^*(A\cap L)-{V(A)\over \zeta
(d)}\Big)^{2} d \nu (L)\le\alpha V(A).$$ Here $\mathscr L(k)=\{
L\in\mathscr L:d(L)\le k\}$, $\zeta(\cdot)$ denotes the Riemann
zeta-function, and $\alpha
>0$ is a constant depending on $k$ and $d$.
\end{form}

The main step of the proof is to show the following proposition:

\begin{form}
\item[\hbox to 1.1cm{\rm (4)\hfill}]
Let $k=1,2,\dots$ Then for $\nu$-almost every lattice
$L\in\mathscr L(k)$ the set $B$ contains infinitely many pairwise
disjoint $d$-tuples of linearly independent points of $L$.
\end{form}
To prove this, let $0=\varrho_{0}<\varrho_{1}<\dots $ be such that
$$V(B_{n})>2^{d}\zeta (d)n, \text{ where } B_{n}=\{ x\in
B:\varrho_{n-1}<\Vert x\Vert \le \varrho_{n}\}.$$ By (2),

\begin{form}
\item[\hbox to 1.1cm{\rm (5)\hfill}]
for $n=1,2,\dots ,$ there are $2^{d}$ pairwise disjoint Borel
sets $B_{ni}, i=1,\dots , 2^{d},$ of $B_{n}$ such that
$V(B_{ni})=\zeta (d)n$ and no $(d-1)$-dimensional subspaces of
$\mathbb E^{d}$ meets each set $B_{ni}$.
\end{form}
Consequently (3) implies that

\begin{form}
\item[\hbox to 1.1cm{\rm (6)\hfill}]
$\displaystyle \int\limits_{\mathscr L(k)} (\#^*(B_{ni}\cap
L)-n)^{2} d\nu (L)\le\alpha\,\zeta (d)n.$
\end{form}
By (5) and (3) the sets $\mathscr L_{ni}=\{ L\in\mathscr
L(k):\#^*(B_{ni}\cap L)=0\},\, i=1,\dots , 2^{d},$ are Borel. It
thus follows from (6) that $n^{2} \nu (\mathscr
L_{ni})\le\alpha\,\zeta (d)n$, or

\begin{form}
\item[\hbox to 1.1cm{\rm (7)\hfill}]
$\displaystyle \nu (\mathscr L_{ni})\le {\alpha\,\zeta (d)\over
n}.$
\end{form}
The set $\mathscr L_{n}=\mathscr L_{n1}\cup\cdots\cup\mathscr
L_{n2^{d}}$ is Borel and consists of all lattices $L\in\mathscr
L(k)$ such that at least one of the sets $B_{ni}$ contains no
primitive point of $L$. Hence $\mathscr L(k)\backslash\mathscr
L_{n}$ is the set of all lattices $L\in\mathscr L(k)$ such that
each set $B_{ni}$ contains a primitive point of $L$. Hence (5)
shows that

\begin{form}
\item[\hbox to 1.1cm{\rm (8)\hfill}]
for any lattice $L\in\mathscr L(k)\backslash\mathscr L_{n}$, the
set $B_{n}$ contains a $d$-tuple of linearly independent points
of $L$.
\end{form}
By (7),

\begin{form}
\item[\hbox to 1.1cm{\rm (9)\hfill}]
$\displaystyle \nu (\mathscr L_{n})\le {\alpha \,2^{d} \zeta
(d)\over n}.$
\end{form}
By definition the sets $B_{n}$ are pairwise disjoint subsets of
$B$. Hence (8) implies that

$$\begin{array}{ll}
\{ L\in\mathscr L(k):&B \text{ contains infinitely many pairwise
disjoint}\,\, d\text{-tuples }\\
&\text{of linearly independent primitive points of } L\}\\
\supset \{ L\in\mathscr L(k):&\text{for infinitely many } n,
\text{ the set } B_{n} \text{ contains a } d\text{-tuple }\\
&\text{of
linearly independent primitive points of } L\}\\
\supset\{ L\in\mathscr L(k):&\text{for infinitely many } n \text{
the lattice } L \text{ is not contained in }\mathscr L_{n}\}\\
\multicolumn{2}{l}{\displaystyle =\bigcap\limits_{m=1}^{\infty}
\bigcup\limits_{n=m}^{\infty} (\mathscr L(k)\backslash\mathscr
L_{n})=\mathscr L(k)\backslash \bigcup\limits_{m=1}^{\infty}
\bigcap\limits_{n=m}^{\infty} \mathscr L_{n}.}
\end{array}$$
Since by (9),
$$\nu \Big( \bigcap\limits_{n=m}^{\infty} \mathscr L_{n}\Big)=0
\text{ for } m=1,2,\dots , \text{ and thus } \nu \Big(
\bigcup\limits_{m=1}^{\infty} \bigcap\limits_{n=m}^{\infty}
\mathscr L_{n}\Big)=0,$$ the proof of (4) is complete.

Since $\mathscr L= \bigcup\mathscr L(k)$, Theorem~1 for $d\ge 3$
is an immediate consequence of (4).

Assume now that $d=2$ and let $\xi$ be the measure on $\mathscr
L$ used by Schmidt \cite{Schmidt60}, p. 525. A result of Schmidt
\cite{Schmidt60}, p. 526/7, shows that (3) continues to hold, but
with the following weaker inequality:
$$\int\limits_{\mathscr L(k)} \Big(\#^*(A\cap L)-{V(A)\over \zeta
(2)}\Big)^{2} d\xi (L)\le\beta\, V(A)\log_{2}V(A).$$ Using this,
we see that in the case $d=2$ the proof is, in essence, the same
as the above proof for $d\ge 3$. Finally, note that the sets of
measure $0$ with respect to $\xi$ and $\nu$ coincide.$\Box$

\section{Proof of the Lemma}

Let $f_{n}$ and $f$ be the distance functions of $S_{n}$ and $S$,
respectively. Since distance functions are positively homogeneous
of degree $1$, and $f_{n}\to f$ uniformly for $\Vert x\Vert\le
1$, we have that $f_n\to f$ uniformly on each bounded set in
$\mathbb E^{d}$. This yields the following statement:

\begin{form}
\item[\hbox to 1.1cm{\rm (10)\hfill}]
Let $l_{n}, l\in\mathbb E^{d}$ be such that $l_{n}\to l$. Then
$f_{n}(l_{n})\to f(l).$
\end{form}
The following claims are simple consequences of the convergence
$L_{n}\to L$, see \cite{GruberLekkerkerker87}, p. 178/9:

\begin{form}
\item[\hbox to 1.1cm{\rm (11)\hfill}]
Given $l\in L$, there are $l_{n}\in L$ such that $l_{n}\to l$.
\end{form}

\begin{form}
\item[\hbox to 1.1cm{\rm (12)\hfill}]
If $l_{n}\in L_{n}$ and $l\in\mathbb E^{d}$ such that $l_{n}\to
l$, then $l\in L$.
\end{form}

(i): Let $\varepsilon >0$. By the definition of successive minima
one can show that there are linearly independent lattice points
$l_{1},\dots , l_{d}\in L$ such that

\begin{form}
\item[\hbox to 1.1cm{\rm (13)\hfill}]
$\max\{ f(l_{1},\dots , f(l_{i})\}\le\lambda_{i}(S,L)+\varepsilon
.$
\end{form}
By (11) we may choose points $l_{nj}\in L_{n}, j=1,\dots , d,$
such that

\begin{form}
\item[\hbox to 1.1cm{\rm (14)\hfill}]
$l_{nj}\to l_{j}.$
\end{form}
Since $l_{1},\dots , l_{d}$ are linearly independent, it follows
that

\begin{form}
\item[\hbox to 1.1cm{\rm (15)\hfill}]
$l_{n1},\dots , l_{nd}\in L_{n}$ are also linearly independent for
all sufficiently large $n$.
\end{form}
Hence, by the definition of $\lambda_{i}$ together with (15),
(14), (10), and (13) we obtain that
\begin{align*}
\lambda_{i}(f_{n},L_{n})&\le\max\{ f_{n}(l_{n1}),\dots ,
f_{n}(l_{ni})\}\le\max\{ f(l_{1},\dots , f(l_{ni})\}+\varepsilon\\
&\le\lambda_{i}(f,L)+2\varepsilon \text{ for all sufficiently
large } n,
\end{align*}
Since $\varepsilon >0$ was arbitrary, this concludes the proof of
claim (i).

(ii): Let $0<\varepsilon <1$. Since $f_{n}\to f$ uniformly on $\{
x:\Vert x\Vert =1\}$ and $f(x)>0$ for $\Vert x\Vert =1$ by the
boundedness of $S$, and $f_{n}, f$ all are continuous and
positively homogeneous of degree 1, there is a constant $\alpha
>0$ such that

\begin{form}
\item[\hbox to 1.1cm{\rm (16)\hfill}]
$\alpha\Vert x\Vert\le (1-\varepsilon ) f(x)\le f_{n}(x)$ for all
$x\in\mathbb E^{d}$ if $n$ is sufficiently large.
\end{form}
For such $n$ we have that $f(x)>0$ for $x\ne o$, hence $S_{n}$ is
bounded. The definition of $\lambda_{i}$ then yields that

\begin{form}
\item[\hbox to 1.1cm{\rm (17)\hfill}]
for all sufficiently large $n$, there are linearly independent
points $l_{n1},\dots , l_{nd}$ in $ L_{n}$, such that
$\lambda_{i}(f_{n},L_{n})=\max\{ f_{n}(l_{n1}),\dots ,
f_{n}(l_{ni})\}$ for $ i=1,\dots , d.$
\end{form}
(16), (17), (1) and (i) together imply that

\begin{form}
\item[\hbox to 1.1cm{\rm (18)\hfill}]
$\displaystyle \Vert l_{ni}\Vert \le {1\over \alpha}
f_{n}(l_{ni})\le {1\over \alpha} \lambda_{i}(f_{n},L_{n})\le
{1\over \alpha} \lambda_{d}(f_{n},L_{n})$

$\displaystyle \phantom{\Vert l_{ni}\Vert}\le {1\over \alpha}
\lambda_{d}(f,L)+\varepsilon$ for all sufficiently large $n$.
\end{form}
For all sufficiently large $n$, the vectors $l_{n1},\dots ,
l_{nd}$ are linearly independent by (17). Consequently, $|\det
(l_{n1},\dots , l_{nd})|$ is an integer multiple of $d(L)$. By
assumption, $L_{n}\to L$. Hence $d(L_{n})\to d(L)$. Combining
this, it follows that

\begin{form}
\item[\hbox to 1.1cm{\rm (19)\hfill}]
$|\det (l_{n1},\dots , l_{nd})|\ge d(L_{n})\ge (1-\varepsilon
)\,d(L)$ for all sufficiently large $n$.
\end{form}
By (18), all the sequences $(l_{n1}),\dots , (l_{nd})$ are
bounded. Fix an index $i=1,\dots , d$. By considering a suitable
subsequence of $1,2,\dots ,$ and re-numbering, if necessary, we
may suppose that

\begin{form}
\item[\hbox to 1.1cm{\rm (20)\hfill}]
$\displaystyle \liminf_{n\to\infty} \lambda_{i}(f_{n},L_{n})$ is
the same as for the original sequence,
\end{form}
and $l_{n1}\to l_{1},\dots , l_{nd}\to l_{d},$ say. By (10), (12)
and (19) the latter implies that

\begin{form}
\item[\hbox to 1.1cm{\hfill}]
$f_{n}(l_{n1})\to f(l_{1}),\dots , f_{n}(l_{nd})\to f(l_{d}),
l_{1},\dots , l_{d}\in L,\\ |\det (l_{1},\dots , l_{d})|\ge
(1-\varepsilon )\, d(L)>0.$
\end{form}
In particular, $l_{1},\dots , l_{d}$ are linearly independent.
Using (17) and the definition of $\lambda_{i}$, it follows that
\begin{align*}
\lambda_{i} (f_{n},L_{n})=\max\{ (f_{n}(l_{n1}),\dots ,
f_{n}(l_{ni})\}\to\max\{ f(l_{1}),\dots ,
f(l_{i})\}\ge\lambda_{i}(f,L).
\end{align*}
This together with (20) and (i) finally yields (ii).$\Box$

\section{Proof of Theorem~2}

$\phantom{iiiii}$(i): Apply Theorem~1 with $B={1\over k}S,
k=1,2,\dots ,$ to see that for $\nu$-almost all lattices
$\varepsilon S$ contains a $d$-tuple of linearly independent
primitive points of $L$ for any $\varepsilon >0$. Hence
$\lambda_{d}(S,L)=0$ for $\nu$-almost all lattices $L$. In
conjunction with (1), this completes the proof of claim (i).

(ii): Let $\mathscr M_{n}=\{ L\in\mathscr L:\lambda_{d}(S,L)\ge
{1\over n}\} , n=1,2,\dots$ Since $\lambda_{d}(S,\,\cdot\, )$ is
upper semi-continuous by the Lemma, $\mathscr M_{n}$ is closed.
If the interior of $\mathscr M_{n}$ is non-empty, then $\nu
(\mathscr M_{n})>0$ by the definitions of $\nu$ and the topology
on $\mathscr L$, in contradiction to (i). Hence $\mathscr M_{n}$
has empty interior. Being closed, $\mathscr M_{n}$ is nowhere
dense in $\mathscr L$. Hence
$$\bigcup\limits_{n=1}^{\infty} \mathscr M_{n}=\{ L\in\mathscr
L:\lambda_{d}(S,L)>0\} \text{ is meager}.$$ Now note (1) to
conclude the proof of claim (ii).$\Box$

\section{Acknowledgement}

For their helpful hints we are obliged to Professors Groemer and
Welzl and the referees. Emo Welzl, in particular, pointed out to
us the existence of the article of Yao and Yao which is reviewed
neither in the Mathematical Reviews nor in the Zentralblatt. The
first author was supported by the Austrian Science Fund (FWF),
project M821-N12.

\hfill\parbox[t]{6truecm}{ Iskander Aliev\hfill\par School of
Mathematics\hfill\par University of Edinburgh\hfill\par James
Clerk Maxwell Building, King's Buildings\hfill\par Mayfield
Road,\hfill\par Edinburgh EH9 3JZ, UK\hfill\par
I.Aliev@ed.ac.uk\hfill\par Peter M. Gruber\hfill\par
 Forschungsgruppe \hfill\par Konvexe und
Diskrete Geometrie \hfill\par Technische Universit\"at
Wien\hfill\par Wiedner Hauptstra\ss e 8--10/1046\hfill\par
A--1040 Vienna, Austria\hfill\par
peter.gruber@tuwien.ac.at\hfill}

\end{document}